\documentclass[letterpaper, 10pt]{article}
  
  \usepackage{geometry}
\usepackage{textcomp}
\usepackage{subfigure}
\usepackage{xcolor}
\usepackage{url}
\usepackage{dsfont} 
\usepackage{acronym}
\usepackage{flushend}
\usepackage{amsmath,graphicx,amssymb,amsfonts}
\usepackage{cite}
\usepackage{hyperref}
\usepackage{tikz}
\usepackage{dsfont} 
\usepackage{amsmath} 
\usepackage{tabularx}
\usepackage{multirow}
\usepackage{enumitem}
\usepackage{cite}

\hypersetup{
    colorlinks=true,
    linkcolor=blue,
    filecolor=magenta,      
    urlcolor=cyan,
    pdfpagemode=FullScreen,
    }
    
\newcommand{\ee}{\,\ensuremath{\mathbf{e}}}

\newcommand{\nm}[0]{(n,m)}

\acrodef{amfm}[AM-FM]{Amplitude-and Frequency-Modulated}
\acrodef{if}[IF]{Instantaneous Frequency}
\acrodef{ia}[IA]{Instantaneous Amplitude}
\acrodef{mcs}[MCS]{Multi Component Signal}
\acrodef{rqf}[RQF]{Reconstruction Quality Factor}
\acrodef{snr}[SNR]{Signal-to-Noise Ratio}
\acrodef{stft}[STFT]{Short-Time Fourier Transform}
\acrodef{cwt}[CWT]{Continuous Wavelet Transform}
\acrodef{tfr}[TFR]{Time-Frequency Representation}
\acrodef{pb}[PB]{Pseudo-Bayesian}
\acrodef{tls}[TLS]{Total Least-Squares}
\acrodef{fri}[FRI]{Finite Rate of Innovation}
\acrodef{svd}[SVD]{Singular Value Decomposition}
\acrodef{rd}[RD]{Ridge Detector}
\acrodef{sstls}[STLS]{Synchrosqueezed TLS}
\acrodef{vsst}[VSST]{Vertical Synchrosqueezing Transform}
\acrodef{tf}[TF]{time-frequency}
\acrodef{dp}[DP]{Dirac pulses} 
\acrodef{svd}[SVD]{Singular Values Decomposition} 
\acrodef{ssa}[SSA]{Singular Spectrum Analysis}

\usepackage{url}

\usepackage{graphicx}

\begin{document}

\title{Time-Frequency Ridge Estimation of Multi-Component Signals using Sparse Modeling of Signal Innovation}

\author{Quentin Legros, Dominique Fourer.
\thanks{This research was supported by the French ANR ASCETE project (ANR-19-CE48-0001).}
}


\maketitle

\begin{abstract}
This paper presents a novel approach for estimating the modes of an observed non-stationary mixture signal.
A link is first established between the short-time Fourier transform and the sparse sampling theory, where the observations are modeled as a stream of pulses filtered by a known function.
As the signal to retrieve has a finite rate of innovation (FRI), an adapted reconstruction approach is used to estimate the signal modes in the presence of noise.
We compare our results with state-of-the-art methods and validate our approach by highlighting an improvement of the estimation performance in different scenarios.
Our approach paves the way of future FRI-based mode disentangling algorithms.
\end{abstract}

%
%

\maketitle

\section{Introduction}

Complex signals generated by a wide range of physical systems are usually modeled as \ac{mcs}, meaning as a sum of amplitude- and frequency-modulated (AM–FM) sines.
Depending on the application, it is often necessary to extract the components for disentangling the signal modes.
To this end, the observed signals are often projected in a \ac{tf} plane \cite{fourer2017c}. Indeed, linear \ac{tf} transforms, such as the \ac{stft}, have received an increasing interest over the last decades for their computational tractability and their ability to highlight the signals components.
The \ac{tfr} of a \ac{mcs} offers an elegant framework where the \ac{if} trajectory of each mode can be observed as a ridge in a 2D plane.
These ridges gather relevant information to denoise the signals, extract specific information, or to separate the sources \cite{meignen2012new,fourer2018b}.
Although estimating the modes of a \ac{mcs} is simple, the presence of noise can avoid the use of classical methods \cite{flandrin2004empirical}.
The \ac{if} estimation performance depends on the presence of outliers and on the quality of the observations \cite{stankovic2001}.
Techniques have been proposed to estimate the \ac{if} by interpolating local maxima along the frequency axis \cite{carmona1997characterization,carmona1999multiridge}.
More recently in \cite{laurent2020novel}, authors extract the regions corresponding to relevant ridge portions from the spectrogram to interpolate remaining pieces of ridges split by the presence of external noise.
In \cite{legros:hal-03259013}, a \ac{pb} framework was used to sequentially estimates the \ac{if} at each time instant using an assumed density filtering approach.
Estimating the ridges position can be achieved using the whole spectrogram, without resorting to sequential ridge tracking \cite{brevdo2011synchrosqueezing,flandrin2015time}.
In \cite{polisano2019}, the problem of lines recovery from degraded images is addressed using a Prony approach on an enhanced version of the \ac{tfr}. This method is well adapted to the \ac{mcs} modes estimation problem, even though deblurring the \ac{tfr} could destruct the ridges with important frequency variation. 
In this work, we introduce a novel approach for estimating the \ac{if} of a \ac{mcs} from its spectrogram. A link is first established with the sparse sampling \cite{VeMaBl02,LIDAR_FRI} theory through modeling of the ridges position. The observed spectrogram can then be viewed as an altered version of a stream of \ac{dp}, filtered by a function that depends on the \ac{tfr} analysis window. The signal of interest thus has a \ac{fri} \cite{VeMaBl02,LIDAR_FRI} whose reconstruction is a known problem.
The presence of spurious noise events is tackled by a \ac{tls} approach \cite{BlDrVe08} instead of the classical Prony method \cite{plonka2019reconstruction}.
One of the advantage of this approach is that the final \ac{if} estimates do not belong to the \ac{tf} resolution grid. Indeed, the reconstruction strategy is able to retrieve real estimates that do not depend on the \ac{tfr} modalities.
We also discuss the use of the \ac{vsst} \cite{oberlin2015second} as an alternative sharpened representation to signals spectrogram.
%
The main contributions of the paper are:
\begin{itemize}
    \item A novel sparse observation model for \ac{mcs} spectrograms.
    \item An \ac{fri} reconstruction strategy, possibly combined with synchrosqueezing.
    \item A method to perform \ac{if} estimation whose performance is independent 
    of the \ac{tfr} resolution.
\end{itemize}
%
This paper is organized as follows. In Section~\ref{sec:problem}, we introduce the problem addressed in this work.
Section~\ref{sec:rec} presents the reconstruction strategy. The performance of the proposed method is comparatively assessed in Section~\ref{sec:results} through numerical experiments.
Conclusions and future work are finally reported in Section~\ref{sec:conclusion}.

\section{Observation model} \label{sec:problem}
Let $x$ be a discrete-time finite-length mixture made of $K$ superimposed \ac{amfm} components expressed as:
\begin{equation} 
 x(n) = \sum_{k=0}^{K-1} x_k(n) \quad \text{, with } x_k(n) = \alpha_k(n)\ee^{2\pi j\phi_k(n)}
\end{equation}
with time instant $n=\{0,1,\dots,N-1\}$, $j^2=-1$ and where $\alpha_k(n)$ and $\phi_k(n)$ are respectively the time-varying amplitude and phase of the $k$-th component.
In the remainder, $K$ is assumed to be known or estimated \cite{Sucic2011}.
The \ac{stft} of signal $x$, using a Gaussian analysis window $\theta_{L}(n)=\frac{1}{\sqrt{2\pi}L}e^{-\frac{n^2}{2L^{2}}}$ with time spread controlled by $L$, can be defined at each time (resp. frequency) instant $n$ (resp. $m\in [0,M-1]$) as:
\begin{equation}\label{eq:stft}
F^{\theta_{L}}_x\nm = \displaystyle\sum_{l=-\infty}^{+\infty} x(l) \theta_{L}(n-l)^* e^{-j\frac{2\pi l m}{M}}
\end{equation}
with $z^*$ the complex conjugate of $z$. Let the square modulus of the \ac{stft} $|F_x^{\theta_{L}}\nm|^2$ be the spectrogram of $x$. We denote this spectrogram by a matrix $\boldsymbol{S}$, whose columns are denoted by $\boldsymbol{s}_{n}$ for a fixed time instant $n$, and $s_{n,m}$ is the spectrogram evaluated at the time $n$ and frequency $m$.
In this work, we are interested in estimating the ridge positions associated with the \ac{if} of the components from a \ac{tfr}.
Note that in the absence of external spurious noise, the spectrogram column $\boldsymbol{s}_{n}$ is assumed to be known if the \ac{if}s $\phi'_{k}(n)$ and \ac{ia} $\alpha_{k}(n)$ are known for each component $k$ at time $n$.
The signal spectrogram can be modeled as:
\begin{equation} \label{eq:Model}
    s_{n,m} \approx \sum\limits_{k=0}^{K-1} a_{k}(n) g(m-\phi'_{k}(n))
\end{equation}
where $a_{k}(n) = \alpha_{k}^2(n)$ and $g(m)=\ee^{-\left(\frac{2\pi m L}{M}\right)^{2}}$ is the squared modulus of the Fourier transform of $\theta_{L}$.
The model in Eq.~\eqref{eq:Model} is only valid for Gaussian analysis window and neglects the components modulation rate. The signal $\boldsymbol{s}_{n}$ can thus be modeled as a mixture of Gaussian functions (since $\theta_{L}$ is Gaussian) whose mean (resp. weight) is controlled by the components \ac{if} (resp. \ac{ia}).
In this work, we address the estimation problem of 1D sparse signals, performed independently $\forall n\in [0,N-1]$, whose only non-zero values indicate the ridges positions and amplitudes of each ridge at time $n$.
The signal to be retrieved can thus be expressed for a fixed time instant $n$ as the following stream of \ac{dp}:
\begin{equation} \label{eq:SoD}
    f_{n}(m) = \sum\limits_{k=0}^{K-1} a_{k}(n) \delta(m-\phi'_{k}(n))\quad.
\end{equation}
The observation modeled in Eq.~\eqref{eq:Model} is then the stream of \ac{dp} in Eq.~\eqref{eq:SoD} filtered by $g$.
We thus address in this work the estimation problem of the \ac{dp}s positions from the observations $\boldsymbol{s}_{n}$. While we discuss \ac{ia} estimation in Section \ref{sec:rec}, this problem remains out of the scope of this work.
Streams of \ac{dp} are signals known to have a \ac{fri} \cite{VeMaBl02,BlDrVe08}.
The models in Eqs.~\eqref{eq:Model} and \eqref{eq:SoD} are limited in the presence of closely spaced components involving cross-terms and oscillations. Nevertheless, this problem can be circumvented by modifying the \ac{tf} resolution, or by resorting to adapted existing approaches \cite{MEIGNEN2019268}. 


\section{Reconstruction} \label{sec:rec}

\underline{Noiseless case: }The restoration of sparse signals in that context has already been studied over the past few years \cite{VeMaBl02,BlDrVe08,hua1990matrix,MaZh93}.
The location and weights of the $K$ \ac{dp}s can indeed be retrieved knowing only $2K+1$ Fourier series coefficients of $f_{n}$ \cite{VeMaBl02}.
From Eq.~\eqref{eq:Model} we obtain:
%
%
\begin{equation} \label{eq:coeff}
    \begin{split}
        s_{n,m} &= \sum\limits_{k=0}^{K-1} a_{k}(n) \sum\limits_{\lambda=-\infty}^{\infty} \hat{g}(\lambda) \ee^{ \frac{j 2\pi \lambda  (m-\phi'_{k}(n))}{M}}\\
        &= \sum\limits_{\lambda=-\infty}^{\infty} \hat{g}(\lambda) \underbrace{\sum\limits_{k=0}^{K-1} a_{k}(n) \ee^{ \frac{-j 2\pi \lambda  \phi'_{k}(n)}{M}}}_{\hat{f}_{n}(\lambda)} \ee^{ \frac{j 2\pi \lambda  m}{M}}
    \end{split}
\end{equation} 
with $\hat{g}$ (resp. $\hat{f}_{n}$) the discrete Fourier transform of $g$ (resp. $f_{n}$).
In order to avoid the use of an infinite sum, a bandlimited approximation can be used in Eq.~\eqref{eq:coeff} such that only $2M_{0}+1$ Fourier series coefficients
are kept \cite{LIDAR_FRI}:
%
\begin{equation} \label{eq:approx}
        s_{n,m} \approx \sum\limits_{\lambda=-M_{0}}^{M_{0}} \hat{g}(\lambda) \hat{f}_{n}(\lambda) \ee^{ \frac{j 2\pi \lambda  m}{M}}
\end{equation}
which rewrites matrix-wise (with the Fourier series coefficients $\hat{f}_{n}(\lambda), \lambda\in[-M_{0},M_{0}]$ of $f_{n}$) as:
\begin{equation} \label{least_square_app}
    \boldsymbol{s}_{n} = \boldsymbol{V} \boldsymbol{D}_{g} \hat{f}_{n} \Leftrightarrow \hat{f}_{n} = \boldsymbol{D}_{g}^{-1} \boldsymbol{V}^{-1} \boldsymbol{s}_{n}
\end{equation}
where $[\boldsymbol{V}]_{m,\lambda} = e^{j2\pi(\frac{m\lambda}{M})}$ is a $(M \times 2M_{0}+1)$ matrix and $\boldsymbol{D}_{g}$ is a diagonal matrix gathering the discrete time Fourier series coefficients of $g$ in $[-M_{0},M_{0}]$. Note that $\boldsymbol{V}$ is an invertible Vandermonde matrix if $M = 2M_{0}+1$.

The second part of the reconstruction method is devoted to retrieve the locations of the \ac{dp} from $\hat{f}_{n}$.
Once $\hat{f}_{n}$ is computed, the annihilating filter method is used to recover the locations $\phi'_{k}(n)$. Let $\boldsymbol{h}$ be a filter that annihilates $\hat{f}_{n}$ such that:
%
\begin{equation} \label{eq:convolutionVandermonde}
    \begin{split}
       (\hat{f}_{n} \ast \boldsymbol{h})(l) &= \sum_{i \in \mathds{Z}} h(i) \hat{f}_{n}(l-i)=0 \\
        &=\sum_{i \in \mathds{Z}} h(i) \sum\limits_{k=0}^{K-1} a_{k}(n) \ee^{ \frac{-j 2\pi (l-i)  \phi'_{k}(n)}{M}}\\
        &=\sum_{k=0}^{K-1} a_{k}(n) \ee^{ \frac{-j 2\pi l \phi'_{k}(n)}{M}} \underbrace{\sum_{i \in \mathds{Z}} h(i)  \ee^{ \frac{j 2\pi i  \phi'_{k}(n)}{M}}}_{H\left(e^{\frac{-j 2\pi \phi'_{k}(n)}{M} }\right)} = 0
    \end{split}
\end{equation}
with $H(z)$ the Z-transform of $\boldsymbol{h}$, whose roots are $e^{\frac{-j 2\pi \phi'_{k}(n)}{M} }$.
Knowing $\boldsymbol{h}$, the $\phi'_{k}(n)$ can thus be retrieved through $H$ roots.
To compute $\boldsymbol{h}$ we assume $h(0)=1$ according to \cite{5686950,maravic2004channel,baechler2017sampling,pan2016towards}, allowing the linear system in Eq.~\eqref{eq:convolutionVandermonde} to be rewritten as follows such that the Toeplitz matrix $\boldsymbol{A}$ is of rank $K$.
%
\begin{equation} \label{eq:Yule-walker_system}
    \underbrace{\!\begin{pmatrix}
        \!\hat{f}_{n}(0) & \!\cdots & \!\hat{f}_{n}(-K+1) \\
        \!\hat{f}_{n}(1) & \!\cdots & \!\hat{f}_{n}(-K+2) \\
         \!\vdots & \!\ddots &  \!\vdots \\
        \!\hat{f}_{n}(K-1) & \!\cdots & \!\hat{f}_{n}(0)
    \end{pmatrix}}_{\boldsymbol{A}}
    \!\begin{pmatrix}
        \!h(1) \\
        \!h(2)\\
        \!\vdots \\
        \!h(K)
    \end{pmatrix}
    \!= \!-\!
    \!\begin{pmatrix}
        \!\hat{f}_{n}(1) \\
        \!\hat{f}_{n}(2)\\
        \!\vdots \\
        \!\hat{f}_{n}(K)
    \end{pmatrix}
\end{equation}
%
This system of Yule-Walker equations has a unique solution if at least $2K$ values of $\hat{f}_{n}$ are known.
Although this is not the focus of our work, this reconstruction approach allows for \ac{ia} estimation. Indeed, from the definition of $\hat{f}_{n}$ in Eq.~\eqref{eq:coeff} once the \ac{if} are estimated we have
%
\vspace{0pt}
\begin{equation}\label{eq:Vandermonde_weight}
    \!\begin{pmatrix}
        \!W_{0,0} & \!\cdots & \!W_{0,K-1} \\
         \!\vdots & \!\ddots & \!\vdots \\
        \!W_{K-1,0} & \!\cdots & \!W_{K-1,K-1} 
    \end{pmatrix}
    \!\begin{pmatrix}
        \!a_{0}(n) \\
        \!a_{1}(n)\\
        \!\vdots \\
        \!a_{K-1}(n) 
    \end{pmatrix}
    \!=\!
    \!\begin{pmatrix}
        \!\hat{f}_{n}(0) \\
        \!\hat{f}_{n}(1)\\
        \!\vdots \\
        \!\hat{f}_{n}(K-1) 
    \end{pmatrix}
\end{equation}
with $W_{p,q}=e^{\frac{-j2\pi p \phi'_{q}(n)}{M}}$. This remains to solve again a Vandermonde system having a unique solution.
%
The presented reconstruction differs from \cite{VeMaBl02}, since $g$ cannot be reduced to a cardinal sine function, but is more similar to the methodology introduced in \cite{LIDAR_FRI} whose development is simpler.
Note that the proposed strategy is not limited to the use of a Gaussian window, and thus any kernel satisfying the Strang-Fix condition \cite{StrangFix} can be used in Eq.~\eqref{eq:stft}.
Moreover, the final estimates are real valued, despite the method assumes that only samples of the observations are known. 

\noindent
\underline{Presence of noise: }
The approach \cite{VeMaBl02} provides, under some assumptions, perfect reconstruction of signals from their sampled and filtered version. Nevertheless, the presence of noise often results in the collapsing of the signal restoration.
The problem of estimating the \ac{if} of signal components requires a particular attention to the presence of noise.
The main limitation when dealing with the presence of outliers is due to the inversion of the matrix $\boldsymbol{V}$ in Eq.~\eqref{least_square_app}, which is sensitive to model mismatch.
Inverting the Vandermonde matrix in Eq.~\eqref{least_square_app} thus becomes an ill-posed problem, avoiding a correct use of the Prony method in Eq.~\eqref{eq:convolutionVandermonde}, since small perturbations in the measurements can result in important reconstruction errors. Different methods have been proposed to circumvent this problem. From them, the \ac{tls} algorithm \cite{BlDrVe08} is known for its simplicity and low computational complexity.
The \ac{tls} method works as an alternative to the Prony approach, in order to compute an approximate annihilating filter.
While this filter will not exactly annihilate $\hat{f}_{n}$, it is a minimizer of $\|\boldsymbol{A}\boldsymbol{h}\|^{2}$, constrained by $\|\boldsymbol{h}\|^{2}=1$.
This minimization of $\|\boldsymbol{A}\boldsymbol{h}\|^{2}$ is a known problem that can be performed by computing the \ac{svd} of $\boldsymbol{A}$ \cite{Cadzow1988SignalEC}, and by setting $\boldsymbol{h}$ to the eigenvector associated with the smallest eigenvalue.
Since model mismatch can significantly limit the performance of the reconstruction method due to the use of a single filtering kernel $g$, we propose to resort to the \ac{vsst} \cite{oberlin2015second} that reallocates the ridges energy to the \ac{if}. 
This method sharpens the \ac{tfr} by concentrating the energy belonging to each ridge on its component \ac{if}. 
Such a representation allows a single filtering kernel $g$ to be used for estimating all signal components. The proposed method is denoted FRI SST when applied on the signal \ac{vsst} instead of its spectrogram.


\section{Results} \label{sec:results}
\begin{figure}[ht!]
\centering
\subfigure{\includegraphics[trim={48pt 0pt 72pt 18pt}, clip,width=0.35\textwidth]{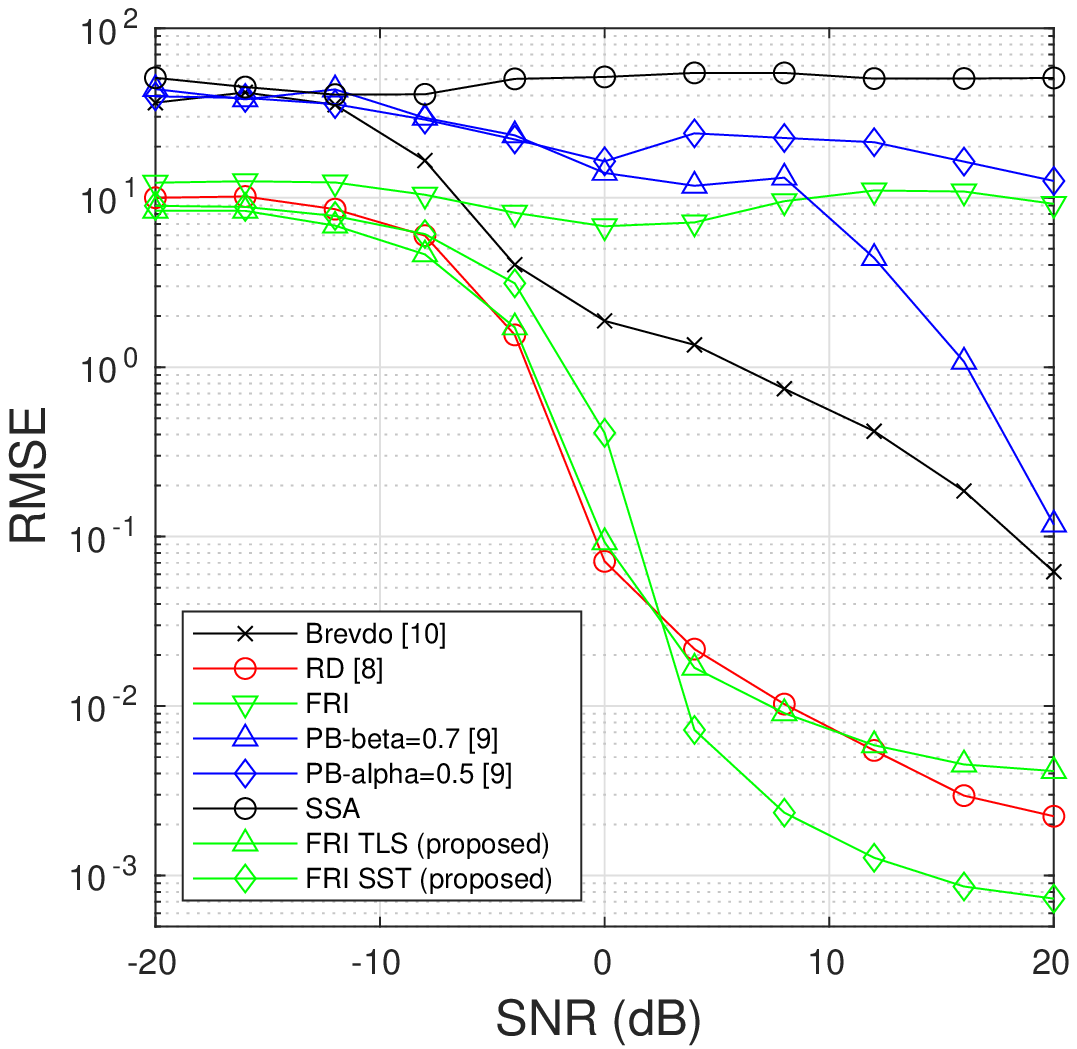}\label{fig:Comp_RMSE}}
\subfigure{\includegraphics[trim={48pt 0pt 72pt 18pt}, clip,width=0.35\textwidth]{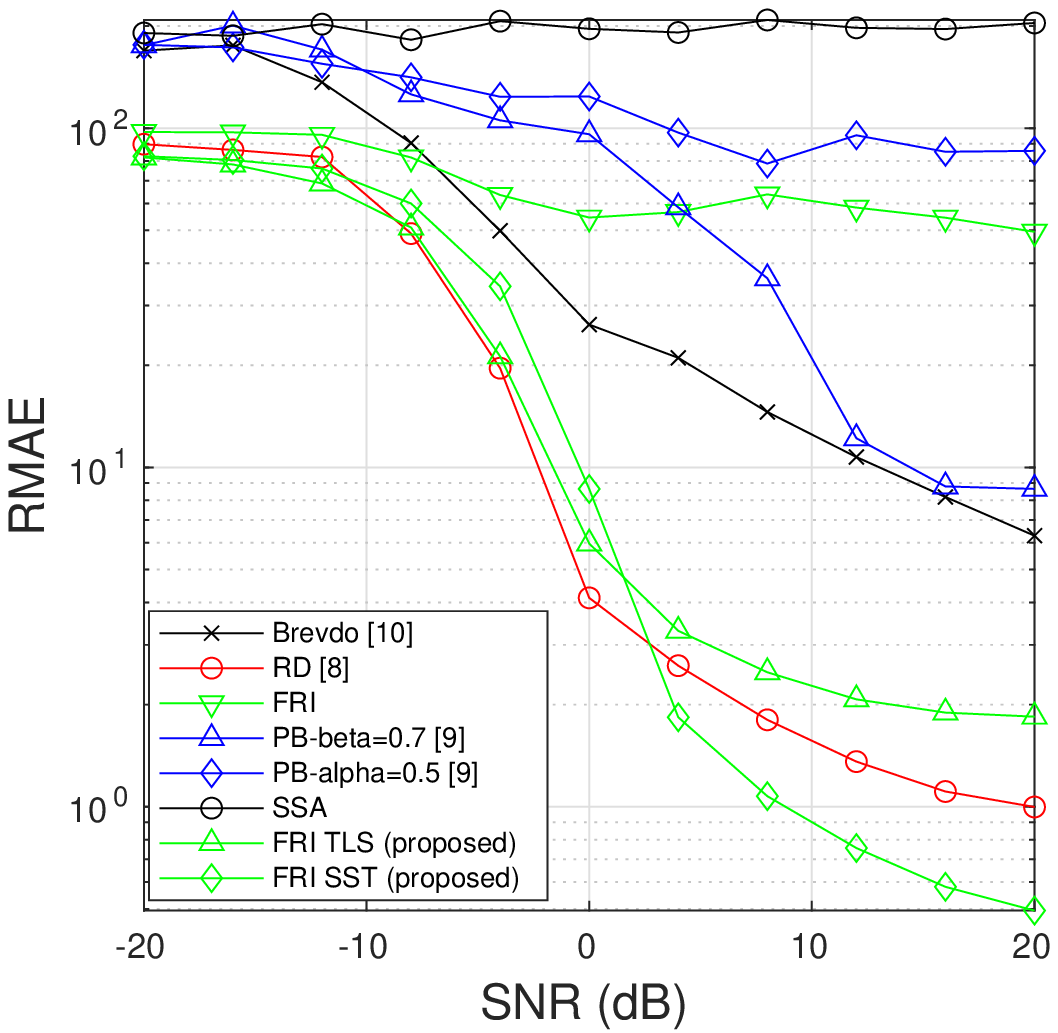}\label{fig:Comp_RMAE}}
\caption{RMSE (top) and RMAE (bottom) of the normalized \ac{if} (averaged over 100 realizations of noise) obtained with the different competing methods for a varying \ac{snr}.}
\end{figure}

%
%
\indent \textbf{\large{1)}  } We assess the \ac{if} estimation performance of the proposed approach\footnote{Codes freely available at \href{https://codeocean.com/capsule/7022037/tree/v1}{Code Ocean} (DOI:10.24433/CO.6654871.v1).}, called FRI TLS, applied on a \ac{mcs} made of three modes: a sinusoid, a linear chirp and a sinusoidally-FM component. 
For the experiments, we compute the \ac{stft} using the ASTRES toolbox \cite{fourer2017c}, with $N=500$ and $L=20$. 
Moreover, we control the \ac{snr} by adding a white Gaussian noise to the \ac{mcs} to simulate the presence of a spurious signal.
We compare the performance of our proposed methods FRI TLS and FRI SST with state-of-the-art approaches: Brevdo \cite{brevdo2011synchrosqueezing}, Noiseless case \ac{fri} presented in Section \ref{sec:rec} (denoted FRI), \ac{rd} \cite{laurent2020novel} and \ac{pb} \cite{legros:hal-03259013} methods. 
We manually set $g$ to a Gaussian function with standard deviation of $0.5$ when resorting to FRI SST, which empirically provides the best correlation with the data.
The estimation performance of the ridges positions is assessed both using the relative mean squared error $\text{RMSE} = \sum\limits_{k=0}^{K-1} \left[\sum\limits_{n=0}^{N-1} \frac{\left(\bar{m}_{n,k}-\hat{m}_{n,k}\right)^{2}}{M^{2}}\right]$, and the relative mean absolute error $\text{RMAE} = \frac{1}{M} \sum\limits_{k=0}^{K-1} \sum\limits_{n=0}^{N-1} \left|\bar{m}_{n,k}-\hat{m}_{n,k}\right|$,
where $\bar{m}_{n,k}$ (resp. $\hat{m}_{n,k}$) is the actual (resp. estimated) normalized \ac{if} of the $k$-th component in the $n$-th time instant.
The RMAE (resp. RMSE) aims to penalize small (resp. important) errors.
The RMSEs (resp RMAEs) obtained with the competing methods are displayed in Fig.~\ref{fig:Comp_RMSE} (resp in Fig.~\ref{fig:Comp_RMAE}) for a varying \ac{snr}.
%
%
%
From Fig.~\ref{fig:Comp_RMSE}, we observe that the proposed FRI TLS and FRI SST approaches, as well as the \ac{rd} provide the best \ac{if} estimation performance. It can be noticed that the \ac{fri} approach is not able to efficiently recover the position of the ridges even at high SNRs.
The proposed FRI SST performs slightly better than the \ac{rd} for all \ac{snr}, except for \ac{snr}$=0$dB, even though both methods performance are similar.
Note that FRI TLS is robust to variations occurring in the assumed filtering kernel $g$, mainly due to the presence of FM components broadening the observed ridges along the frequency axis \cite{colominas2019}.
Moreover, applying the proposed method on a signal \ac{vsst} produces similar performance than when applied on the spectrogram for \ac{snr}$<0$dB.
The results displayed in Fig.~\ref{fig:Comp_RMAE} aims to highlight the behavior of the competing method in a more local perspective. 
From Fig.~\ref{fig:Comp_RMAE}, a similar behavior to the results in Fig.~\ref{fig:Comp_RMSE} can be observed, where the  FRI TLS, FRI SST and the \ac{rd} obtain the best performances. While their RMAE converge to a close value in the low \ac{snr} regime, the FRI TLS and FRI SST reach the best performance for \ac{snr}$<10$dB, while the \ac{rd} method obtains the lowest error for a \ac{snr}$=0$dB. These results, compared to those displayed in Fig.~\ref{fig:Comp_RMSE}, show that the FRI TLS and FRI SST estimates are close but oscillate around the true \ac{if}. 
No spatial constraint is used to regularize the estimation process conversely to the compared methods, while it would avoid small oscillations and improve the overall performance of the approach.

%
\indent \textbf{\large{2)}  } We then demonstrate the ability of the proposed method to approximate the signal modes by reconstructing separately the $K=3$ components C1, C2 and C3, by filtering the signal \ac{tfr} with a binary filter of width $2\times 10 +1$, centered around the components \ac{if}. The inverse \ac{stft} formula is then applied to each extracted mode whose reconstruction performance is assessed using the \ac{rqf}:
$10 \log_{10}\left( \frac{||x||^2}{||x-\hat{x}||^2}\right)$,
where $x$ (resp. $\hat{x}$) stands for the reference (resp. estimated) signal.
\begin{table}[htb]
\caption{component-wise \ac{rqf} (averaged over 100 realizations of noise) of the competing approaches for different \ac{snr}. The standard deviation of the estimators are displayed in second rows .\label{table:Comp_IF_MCS}}
\begin{center}
\subtable[SNR=10dB]{
\begin{tabular}{l|l|l|l|l|}
\cline{2-5}
& C1 & C2 & C3 & Average \\ 
\hline
\multicolumn{1}{|l|}{RD \cite{laurent2020novel}}    & 
\begin{tabular}[c]{@{}l@{}}$16.73$\\ $\pm 17.42$\end{tabular} &
\begin{tabular}[c]{@{}l@{}}$16.12$\\ $\pm 2.48$\end{tabular} &
\begin{tabular}[c]{@{}l@{}}\textbf{15.68}\\ $\pm 0.54$\end{tabular} &
\begin{tabular}[c]{@{}l@{}}$16.18$\\ $\pm 2.10$\end{tabular} \\ 
\hline
\multicolumn{1}{|l|}{FRI TLS (proposed)}   & 
\begin{tabular}[c]{@{}l@{}}\textbf{20.30}\\ $\pm 0.81$\end{tabular} &
\begin{tabular}[c]{@{}l@{}}$18.84$\\ $\pm 1.15$\end{tabular} &
\begin{tabular}[c]{@{}l@{}}$ 11.40$\\ $\pm 0.6$\end{tabular} &
\begin{tabular}[c]{@{}l@{}}\textbf{16.84}\\ $\pm 0.88$ \end{tabular} \\ 
\hline
\multicolumn{1}{|l|}{FRI SST (proposed)}   & 
\begin{tabular}[c]{@{}l@{}}$20.00$\\ $\pm 0.89$\end{tabular} &
\begin{tabular}[c]{@{}l@{}}\textbf{19.04}\\ $\pm 1.30$\end{tabular} &
\begin{tabular}[c]{@{}l@{}}$ 11.30$\\$\pm 0.53$\end{tabular} &
\begin{tabular}[c]{@{}l@{}}$16.78$\\ $\pm 0.51$ \end{tabular} \\ 
\hline
\end{tabular}}
\subtable[SNR=0dB]{
\begin{tabular}{l|l|l|l|l|}
\cline{2-5}
& C1 & C2 & C3 & Average \\ 
\hline
\multicolumn{1}{|l|}{RD \cite{laurent2020novel}}    & 
\begin{tabular}[c]{@{}l@{}}$9.00$\\ $\pm 1.31$\end{tabular} &
\begin{tabular}[c]{@{}l@{}}$8.21$\\ $\pm 1.02$\end{tabular} &
\begin{tabular}[c]{@{}l@{}}$6.56$\\ $\pm 1.19$\end{tabular} &
\begin{tabular}[c]{@{}l@{}}$7.92$\\ $\pm 1.18$\end{tabular} \\ 
\hline
\multicolumn{1}{|l|}{FRI TLS (proposed)}   & 
\begin{tabular}[c]{@{}l@{}}$8.81$\\ $\pm 1.04$\end{tabular} &
\begin{tabular}[c]{@{}l@{}}\textbf{8.28}\\ $\pm 1.54$\end{tabular} &
\begin{tabular}[c]{@{}l@{}}\textbf{7.37}\\ $\pm 0.99$\end{tabular} &
\begin{tabular}[c]{@{}l@{}}\textbf{8.15}\\ $\pm 1.22$ \end{tabular} \\ 
\hline
\multicolumn{1}{|l|}{FRI SST (proposed)}   & 
\begin{tabular}[c]{@{}l@{}}\textbf{9.16}\\ $\pm 1.63$\end{tabular} &
\begin{tabular}[c]{@{}l@{}}$7.53$\\ $\pm 1.45$\end{tabular} &
\begin{tabular}[c]{@{}l@{}}$6.03$ \\ $\pm 0.93$\end{tabular} &
\begin{tabular}[c]{@{}l@{}}$7.57$ \\ $\pm 1.37$ \end{tabular} \\ 
\hline
\end{tabular}}
\end{center}
\end{table}
%
The \ac{rqf} obtained for \ac{snr}s of $0$dB and $10$dB are displayed in Table~\ref{table:Comp_IF_MCS}, where we discarded all the results except those of the proposed approaches and RD.
The proposed FRI TLS obtains the highest \ac{rqf} in average for the three components for both \ac{snr}s. Moreover, both FRI TLS and FRI SST provide the best performance for all component, except for component C3 with \ac{snr}$=10$dB, due to important modulation rate involving model mismatch between the observed data and $g$. 
For \ac{snr} = $0$dB, the highest \ac{rqf} are obtained using FRI TLS except for the C1. Indeed, the sharpened \ac{tf} representation used with FRI SST circumvents the limitations induced by ridges shape modifications in the presence of noise. This last alternative is particularly effective for sinusoidal components.

\begin{figure}[ht!]
\centering
\includegraphics[trim={48pt 0pt 72pt 19pt}, clip,width=0.35\textwidth]{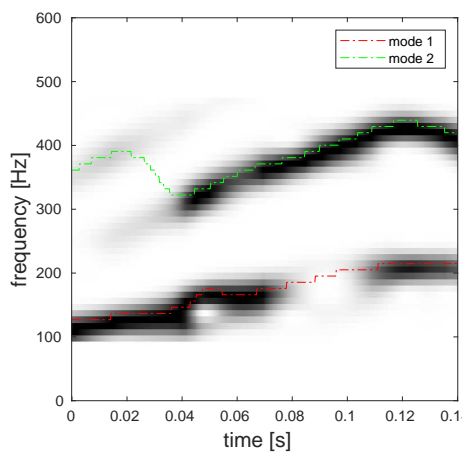}
\caption{Estimation of the first $K=2$ signal components of the speech signal using the proposed \ac{tls} method.\label{fig:speech2}}
\end{figure}

%
\indent \textbf{\large{3)}  } We finally consider a real-world signal made of three components defined on a reduced portion of the time axis in Fig.~\ref{fig:speech2}. Nevertheless, we assume for the presence of only two components to enhance readability of the results and for assessing the \ac{tls} method behavior in that case (\ac{stft} computed using $L=40$).
%
The transition between the two parts of mode 2 is smooth, since the method searches for the filtered \ac{dp} providing the strongest correlation with the data. Minimizing the $\ell_2$-norm as detailed in Section \ref{sec:rec} remains to select the \ac{tf} points minimizing the mean-square error, corresponding to the mean between the two ridges of mode 1.

\section{Conclusion} \label{sec:conclusion}
We have proposed a novel approach for estimating the ridges of a \ac{mcs} in presence of noise.
Our method assumes that the signal to reconstruct is sparse, since for a given time instant, each component is associated with a pulse whose position indicates its \ac{if}.
The spectrogram time slices are viewed as noisy filtered and sampled streams of Dirac pulses, reducing the problem of retrieving \ac{fri} signals.
Since the presence of noise avoids an efficient use of the Prony method, we resort to the total least square  alternative.
The comparison with other state-of-the-art approaches shows a significant improvement of the \ac{if} estimation performance.
A simple alternative involving only an additional transform of the signal was also presented to deal with signal components spreading inconsistently in the frequency plane.
Future work should focus on efficiently disentangling the estimated components and on investigating adapted spatial regularization. 


\end{document}